\documentclass[preprint]{elsarticle}

\usepackage{amssymb,amsmath, color}
\usepackage[mathscr]{eucal}

\usepackage{fancyhdr}
\makeatletter
\def\ps@pprintTitle{%
 \let\@oddhead\@empty
 \let\@evenhead\@empty
 \def\@oddfoot{}%
 \let\@evenfoot\@oddfoot}
\makeatother
\newtheorem{thm}{Theorem}  

\newtheorem{lemma}[thm]{Lemma}

\newproof{pf}{Proof}
\newdefinition{rmk}{Remark}

\newcommand{\sX}{\mathscr{X}}  
\newcommand{\td}[1]{\tilde{#1}}
\newcommand{\R}{\mathbb{R}}

 \newcommand{\sg}{\sigma}
 \newcommand{\Om}{\Omega}

\begin{document}

\begin{frontmatter}



\title{ Twice $Q$-polynomial distance-regular graphs of diameter 4 }

\author[JM]{Jianmin Ma}
\ead{jianminma@yahoo.com}
\author[JK]{Jack H.  Koolen \corref{cor1}}
\ead{koolen@ustc.edu.cn}

\cortext[cor1]{Corresponding author}

\address[JM]{Hebei Key Lab of Computational Mathematics \& Applications  and  College of Math \& Info. Science \\  Hebei Normal University, Shijiazhuang, Hebei 050016, China}

\address[JK]{School of Mathematical Sciences, University of Science and Technology of China,
 Hefei, Anhui 230026, China}

\begin{abstract}
It is known that a distance-regular graph with valency $k$ at least three admits at most two Q-polynomial structures.
In this note we show that all distance-regular graphs with diameter four and valency at least three admitting two $Q$-polynomial structures are either dual bipartite or almost dual imprimitive.
By  the work of Dickie \cite{Dickie} this implies that
any distance-regular graph with diameter $d$ at least  four and valency at least three admitting two $Q$-polynomial structures is, provided it is not a Hadamard graph, 
 either the cube $H(d,2)$ with $d$ even, the half  cube ${1}/{2} H(2d+1,2)$, the folded cube  $\tilde{H}(2d+1,2)$, or 
the dual polar graph on $[^2A_{2d-1}(q)]$ with $q\ge 2$ a prime power.
\end{abstract}

\date{6/17/2013}
\end{frontmatter}

\section{Introduction}
\label{s:intro}

Dickie \cite {Dickie} showed that any distance-regular graph with diameter at least four and valency at least three admits at most two $Q$-polynomial structures,  extending the work of Bannai and Ito who showed it for diameter at least 34.
For brevity we call  a distance-regular graph (or an association scheme) with exactly two $Q$-polynomial structures {\em twice $Q$-polynomial}. 
Furthermore, Dickie \cite {Dickie} classified twice $Q$-polynomial distance-regular graphs with diameter at least five and valency at least three.

\begin{thm} \cite[Theorem 8.1.2]{Dickie}
\label{T:D812} \label{t:Dickie} \label{t:drg2Q}
Let $\varGamma$ be a distance regular graph with diameter $d\ge 5$ and valency $k\ge 3$. Then $\varGamma$ has two $Q$-polynomial structures if and only if $\varGamma$ is one of the following:
\begin{enumerate}[(i)]
\item
the cube $H(d,2)$ with $d$ even;
\item
the half  cube $\frac{1}{2} H(2d+1,2)$;
\item
the folded cube  $\tilde{H}(2d+1,2)$;
\item
the dual polar graph on $[^2A_{2d-1}(q)]$, where $q\ge 2$ is a prime power.
\end{enumerate}
\end{thm}

In this note we show that Theorem 1 can be extended to include the diameter four case. The following result is key to doing so.


\begin{thm} \label{t:main}
Let $\varGamma$ denote a twice $Q$-polynomial distance-regular graph of diameter four and valency at least three. Then  one of the $Q$-polynomial structures  has $a_1^*= a_2^* =a_3^* =0$, that is, this structure is  either dual bipartite or almost dual bipartite.
\end{thm}

As a consequence of Theorems \ref{t:drg2Q}, \ref{t:main} and the  classification of  distance-regular graphs of diameter
at least four that are either  dual bipartite or  almost dual bipartite  \cite{Dickie},  we obtain the following result. 
\begin{thm} \label{mainD=4}
Let $\varGamma$  denote a twice $Q$-polynomial distance-regular graph with diameter $d$ at least $4$ and valency at least $3$. Then $\varGamma$   is one of the following:
\begin{enumerate}[(i)]
\item
the cube H(d,2) with $d$  even;
\item
the half  cube $\frac{1}{2} H(2d+1,2)$;
\item
the folded cube  $\tilde{H}(2d+1,2)$;
\item
the dual polar graph on $[^2 A_{2d-1}(q)]$, where $q\ge 2$ is a prime power;
\item
a Hadamard graph of order $2 \gamma$ with intersection array
$\{2 \gamma, 2 \gamma-1,  \gamma, 1; 1, \gamma, 2\gamma -1, 2 \gamma\}$
 with $\gamma =1$ or  $\gamma$ a positive even integer.

 \end{enumerate}
\end{thm}

\begin{rmk}
It is known that  Hadamard graphs of order $2\gamma$ exist for $\gamma$ any non-negative power of two.
A Hadamard graph of order $2\gamma$ exists if and only if a Hadamard
matrix  of $2\gamma$
 exists \cite[Section 1.8]{BCN}. 
A Hadamard matrix
of order $n$ with $n$ a positive integer is a square  $\{+1,-1\}$ matrix $H$ of order $n$ such that $H H^T = nI$.
 The Hadamard conjecture states that a Hadamard matrix exists if and only if $n =1,2$, or $n$ is a positive integer divisible by $4$.
\end{rmk}

Distance-regular graphs of diameter 2 are  strongly regular graphs, which  possess  two $P$-polynomial and two $Q$-polynomial structures. Any connected distance regular graph with valency two is an ordinary $n$-gon, which can have more than two $Q$-polynomial structures only if $n\ge 7$. 
%
So in the rest of this note, we restrict ourselves to 
 distance-regular graphs with both diameter and valency at least three unless stated otherwise.

\section{Definitions and preliminaries}

In this section, we review some definitions and basic concepts. See the books of Brouwer, Cohen, and Neumaier \cite{BCN} or Bannai and Ito \cite{BI} for more background information.

\subsection{Distance-regular graphs}
Let $\varGamma=(X,R)$ denote a finite, undirected, connected graph,
without loops or multiple edges, with vertex set $X$, edge set
$R$, path-length distance function  $\partial$, and diameter
$d  =  \hbox{max}\left\{ \partial(x,y) \mid x,y\in X\right\}.$
For all $x \in X$ and for all integers $i$, we set
$\varGamma_i(x) =  \lbrace y\in X \mid \partial(x,y)=i\rbrace.$
We abbreviate $\varGamma(x) = \varGamma_1(x)$.
The {\em valency} of a vertex  $x\in X$ is the
cardinality of $\varGamma(x)$.
The graph $\varGamma $ is said to be {\em regular,  valency k}, if
each vertex in  $X$ has valency $k$. Graph $\varGamma$ is said to be
{\em distance-regular} whenever for all integers
$ h, i, j$ $ (0 \le h, i, j\leq d)$, and for all $ x, y \in X $
with $ \partial (x,y) = h$, the number
$
p^h_{i,j}  =  \vert \varGamma_i(x) \cap \varGamma_j(y)\vert
$
is independent of $ x$ and $ y$.
When there is no possibility of confusion,  we write $p^h_{ij}$ instead $p^h_{i,j}$.
The constants $ p^h_{ij} $
are known as the {\em intersection numbers} of $ \varGamma$.

For notational convenience,
set $ c_i= p^i_{1, i-1}$ $ (1 \le i \leq d)$, $ a_i  =  p^i_{1i}$ $
(0 \leq i \leq d)$, $ b_i  =  p^i_{1, i+1}$ $ (0 \le i \le d-1)$, $
k_i  =  p^0_{ii}$ $(0\leq i \leq d)$,  and  define
$ c_0 = 0$, $ b_d = 0$. We note $a_0=0$, $k_1 = b_0$, and   $c_1=1$. We write $k=k_1$. The sequence
$$
\{b_0,b_1,\dots,b_{d-1}; c_1,c_2,\dots,c_d\}
$$
is called the {\em intersection array} of $\varGamma$.

From now on,  $ \varGamma=(X,R)  $ will denote a  distance-regular graph
of diameter $d\geq 3$. Observe that $\varGamma$ is regular with valency
$k$, and that
\begin{equation}
\quad  k =  c_i  + a_i + b_i \qquad (0 \le i \le d).
\label{pre2}
\end{equation}

We now recall the Bose-Mesner algebra.
For each integer $ i$ $ (0 \le i \le d)$, let $A_i$ denote the $i$-th distance matrix
 with
$ x, y $
entry
$$
(A_i)_{xy}= \begin{cases}
1, & \mbox{ if } \partial(x,y) = i,\\
0, & \mbox{ otherwise },
\end{cases}
 \quad (x, y \in X).
$$
Then
\begin{equation}\label{d:scheme}
A_0 = I, \quad A_0 + A_1 + \dots + A_d = J, \quad A_i^t = A_i, \quad A_iA_j = \sum_{k=0}^d p^k_{ij},
\end{equation}
where $J$ is the all-one matrix.

We abbreviate $A = A_1$, and refer to this as the
{\em adjacency matrix} of $\varGamma$.
Let $M$ denote the algebra  generated by  $A$ over the reals $\R$. We refer to $M$ as the
{\em Bose-Mesner algebra} of $\varGamma $. The matrices
 $ A_0, A_1,\ldots, A_d $ form a basis
for $M$. Note that $M$ is closed under the Schur (entry-wise) product $\circ$. So it has a second basis
$ E_0, E_1, \ldots, E_d $ of primitive idempotents which satisfy
\begin{equation}\label{d:idempotent}
E_0 = |X|^{-1} J,\quad
 E_0 + E_1+ \dots E_d = I, \quad
 E^t_i= E_i,\quad
E_i \circ E_j =|X|^{-1}\sum_{k=0}^d q^k_{ij} E_k.
\end{equation}
The numbers $q^k_{ij}$ are nonnegative reals, and are referred to as {\em Krein parameters}. Note the  parameters $q^k_{ij}$ and $p^k_{ij}$ depend on the orderings $(E_i)_i$, and $(A_i)_i$, respectively.

Let  $ \theta_0,\theta_1,\ldots,\theta_d$ denote the real numbers satisfying
$
A = \sum^d_{i=0} \theta_i E_i.
$
We refer to $\theta_i$ as the {\em eigenvalue} of $\varGamma $ associated
with $E_i$, and  call $\theta_0$ the {\em trivial} eigenvalue.
If $\theta_0 > \theta_1>\dots > \theta_d$, then we say $E_0,E_1, \dots, E_d$ is the {\em natural ordering} of the primitive idempotents.

A {\em $d$-class   association scheme} is a pair $(X, \{R_0, \dots, R_d\})$ with $R_i$ symmetric binary relations on $X$ whose adjacency matrices satisfy \eqref{d:scheme}, where $R_i$ has {\em adjacency matrix $A_i$} defined by $(A_i)_{xy}=1$ if $(x,y)\in R_i$ and $0$ otherwise.  It also have a set of primitive idempotents that satisfy \eqref{d:idempotent}.
The property that one of the relations of a $d$-class association scheme forms a distance-regular
graph with diameter $d$ is equivalent to the scheme being {\em P-polynomial}, that is, the
relations $R_0, \dots, R_d$ can be ordered such that every pair of vertices in $R_i$ has distance $i$ in the graph
 $(X, R_1)$ for every $i$. In turn, this is equivalent to the conditions $p^{i+1}_{1i} >0 $ for $0\le i\le d-1$  and
$p^{k}_{1i} =0$ for $k > i+1$.

\subsection{Cosines}
 We now recall the cosines.
Let $\theta $ denote an eigenvalue of $\varGamma $, and  let
$E$ denote the associated primitive idempotent.
Let $\sigma_0, \sigma_1, \ldots, \sigma_d$ denote the real numbers satisfying
\begin{equation}
E = \vert X \vert^{-1}m\sum_{i=0}^d \sigma_iA_i,
\label{pre3}
\end{equation}
where $m$ denotes the multiplicity of $\theta$.
Taking the trace in (\ref{pre3}), we find $\sigma_0=1$.
We call $\sigma_i$  the
{\em i-th cosine} of $\varGamma $ with respect to $\theta $ (or $E$), and call
 $\sigma_0, \sigma_1, \ldots, \sigma_d$
the {\em cosine sequence} of $\varGamma $ associated with $\theta $ (or $E$).

We will need the following basic results.

\begin{lemma}\cite[Section 4.1.B]{BCN} \label{Lem4.3}  \label{cosineseq}            
Let $\varGamma$ denote a distance-regular graph with diameter $ d\geq 3$. Then for any complex numbers
$ \theta,  \sigma_0, \sigma_1, \ldots, \sigma_d$, the following are equivalent.
\begin{enumerate}[(i)]
\item
$\theta $ is an eigenvalue of $ \varGamma $, and	   $ \sigma_0, \sigma_1, \ldots, \sigma_d $ is the	   associated cosine sequence.
\item
$\sigma_0=1$, and
\begin{equation}
	     c_i\sigma_{i-1} +a_i\sigma_i+b_i\sigma_{i+1}
	     = \theta \sigma_i \qquad \qquad (0 \le i \le d),
\label{pre2A}
\end{equation}
             where $ \sigma_{-1}$ and  $\sigma_{d+1} $ are
	     indeterminates.
\item
$\sigma_0=1$, $k \sigma = \theta$, and
            \begin{equation}
             c_i(\sigma_{i-1} - \sigma_i) - b_i(\sigma_i-\sigma_{i+1})
             = k(\sigma-1) \sigma_i \qquad \qquad (1 \le i \le d),
           \label{pre2B}
	   \end{equation}
             where  $\sigma_{d+1} $ is an
	     indeterminate.
\hfill $\Box$
\end{enumerate}
\end{lemma}

The second largest
and minimal eigenvalue of a distance-regular graph turn out to be
of particular interest.
In the next several lemmas, we give some basic information on these
eigenvalues.
\begin{lemma}\cite[Lem. 13.2.1]{Godsil} \label{L5.1} \label{ZL2.6}   
Let $\varGamma$ denote a distance-regular graph with diameter
$ d\geq 3$, and eigenvalues $ \theta_0> \theta_1>\cdots >\theta_d$.
Let $\theta$ denote one of $\theta_1, \theta_d$ and let
$ \sigma_0, \sigma_1, \ldots, \sigma_d$ denote the cosine sequence
for $\theta $.
\begin{enumerate}[(i)]
\item
 Suppose $\theta=\theta_1$. Then $ \sigma_0  > \sigma_1 > \cdots > \sigma_d$.
\item
 Suppose $\theta=\theta_d$. Then for each $i$ $(0 \le i \le d)$,  $ (-1)^i\sigma_i > 0$.
 \hfill $\Box$
\end{enumerate}
\end{lemma}

\begin{lemma} \label{ZL2.0}       
Let $\varGamma$ be  a distance-regular graph with
 diameter $ d\geq 3$ and eigenvalues
$ \theta_0 > \theta_1 > \cdots > \theta_d$.  Then  {\em (i)--(iii)}  below
hold.
\begin{enumerate} [(i)]
\item $0 < \theta_1 < k$.
\item    $a_1-k\le \theta_d <-1$.
\item If $\varGamma $ is not bipartite,  then $ a_1 -k< \theta_d$.
\hfill $\Box$
\end{enumerate}
\end{lemma}

\subsection{$Q$-polynomial property}

Let $\theta_0, \theta_1, \dots, \theta_d$ (or $E_0, E_1, \dots, E_d)$ be a fixed ordering of the eigenvalues (or primitive idempotents) of $\varGamma$. We call this ordering is a {\em $Q$-polynomial structure} if there is a sequence $\sigma =(\sigma_0, \sigma_1, \dots, \sigma_d)$ and polynomial $q_j$ of degree $j$, $j=0,1,\dots,d$, such that
$$
E_j = \sum_{i=0}^d q_j(\sigma_i) A_i;
$$
in this case, $\sigma$ is called a {\em $Q$-sequence} of $\varGamma$ and $E_1$ is called the {\em primary idempotent} for this $Q$-sequence. The graph $\varGamma$ is called {\em $Q$-polynomial} if $\varGamma$ has a $Q$-polynomial structure.

Let $E_0, E_1, \dots, E_d$ be a $Q$-polynomial structure for $\varGamma$. We usually write
$a_i^* = q^i_{1i}$, $b_i^* = q^i_{1, i+1},$, $c_i^*=q^i_{1,i-1}$ and $k_i^* = q^0_{ii}$ for $i=0,1,\dots,d$.

We will need the following theorem of H. Suzuki.

\begin{thm}  \label{t:2Q} \cite{Suzuki98}
 Suppose that $\sX $ is a symmetric  association scheme with a $Q$-polynomial structure
  $E_0, E_1, \dots, E_d$. If $\sX$ is not a polygon and has another   $Q$-polynomial structure, then the new structure is one of the following:
\begin{enumerate}[(I)]
	\item
	$E_0, E_2, E_4, E_6, \dots, E_5, E_3, E_1;$
	\item
	$E_0, E_d, E_1,  E_{d-1},  E_2, E_{d-2}, E_3, E_{d-3}, \dots;$
	\item
	$E_0, E_d, E_2, E_{d-2}, E_4, E_{d-4}, \dots, E_{d-5}, E_5, E_{d-3}, E_3, E_{d-1}, E_1;$
	\item
	$E_0, E_{d-1}, E_2, E_{d-3}, E_4, E_{d-5}, \dots, E_5, E_{d-4}, E_3, E_{d-2}, E_1, E_d;$
	\item
	{$d=5$ and $E_0, E_5, E_3, E_2, E_4, E_1$.}
	\end{enumerate}
Hence, $\sX$ admits at most two $Q$-polynomial structures.
\end{thm}
Case (V) was recently eliminated in  \cite{MaW}.

\subsection{Almost dual primitivity}\label{s:dualprim}

The graph $\varGamma$ is called {\em imprimitive} when some $i$, $1\le i\le d$, the distance-$i$ graph $\varGamma_i = (X,A_i)$ is disconnected. If $\varGamma$ is imprimitive, then by \cite[Theorem 4.2.1]{BCN}, $\varGamma$ is {\em bipartite} (here $\varGamma_2$ is disconnected) or {\em antipodal} (here $\varGamma_d$ is a union of cliques).

A  $Q$-polynomial structure $(E_i)^d_{i=0}$ is called  {\em dual bipartite} if  $a_0^* =a_1^* =\dots = a_d^*= 0$. When there is no possibility of confusion, we also say that graph $\varGamma$ is dual bipartite. Similar comment applies the other concepts to follow immediately.  If $c_i^* = b_{d-i}^*$ for $i=0,1,\dots,d$ and $i\not=\lfloor d/2\rfloor$,  then $\varGamma$ is called {\em dual antipodal}. An imprimitive $Q$-polynomial distance-regular graph is either dual bipartite or dual antipodal (or both).

Now we define the terms of almost dual bipartite/antipodal, introduced by Dickie \cite{Dickie}. A  $Q$-polynomial structure $(E_i)^d_{i=0}$ is called  {\em almost dual bipartite} if $a_0^* =a_1^* =\dots = a_{d-1}^*= 0 \ne a_d^* $; it  is called  {\em almost dual antipodal}  if $q^d_{1d} \ne 0 = q^d_{2d} = \dots = q^d_{dd}$. If $\varGamma$ is almost dual bipartite or antipodal, then it is called {\em almost dual imprimitive}.

For a classification of almost dual imprimitive  distance-regular graphs,  see \cite[Theorem 3.1.4] {Dickie} for the almost dual bipartite case with $d\ge 4$, and \cite[Theorem 2.1.2]{Dickie},\cite{Dickie96} for the dual bipartite case with $d\ge 3$.

\subsection{The tight property}

Now we recall the tight property \cite{Jurisic00}.  A distance-regular graph $\varGamma$ is called {\em tight} if it is not bipartite and the following equality holds
$$
\left(\theta_1 + \frac{k }{a_1 +1}\right) \left(\theta_d + \frac{k}{a_1 +1}\right) = -\frac{k a_1 b_1}{(a_1^2 +1)^2},
$$
where $\theta_1$ and $\theta_d$ are the second largest and the smallest eigenvalues of $\varGamma$, respectively.
\section{Proof of Theorem \ref{mainD=4}}\label{s:proof}

We prove our main theorems in this section. In the rest of this paper, we will fix $\Om$ to be a general distance-regular graph and $\varGamma$ to be a twice $Q$-polynomial distance-regular graph of diameter 4. Let $E_0, \dots, E_4$ and $\td{E}_0,\dots, \td{E}_4$ be $Q$-polynomial structures for $\varGamma$.  The parameters for $(\td{E_i})_i$ will be attached with a tilde.

 We first prove Theorem \ref{t:main}, which is key to the proof of Theorem \ref{mainD=4}.

\subsection{Proof of Theorem \ref{t:main}}

    We first quote some results from Dickie's thesis \cite{Dickie}.
\begin{thm}\cite[Theorem 7.1.1]{Dickie},\cite{Dickie95}
\label{T:D711}
Let $\Om = (X, R)$ denote a distance regular graph with diameter $d\ge 3$. Suppose that $\Om$ admits more than one $Q$-polynomial structure. Then $\Om$ is thin and dual thin.
\hfill $\Box$
\end{thm}

The follow result follows from Theorem \ref{T:D711} and Theorems 4.1.1 and 5.1.1. in \cite{Dickie}.

\begin{thm}\cite[Theorem 5.1.2]{Dickie}
\label{T:D512}
Let $\Om$ be as in Theorem \ref{T:D711} with a $Q$-polynomial structure  $E_0, E_1, \dots, E_d$,  Krein parameters $q^i_{1i}$ and intersection numbers $p^i_{1i}$. Then we have the following implications:
 \begin{align}
 q^{1}_{11} =0 \quad \Rightarrow \quad q^{i}_{1i} =0 \label{eq:dualbip}\\
 q^{1}_{11} \not=0 \quad \Rightarrow \quad q^{i}_{1i} \ne 0 \label{eq:dualNbip}\\
 p^{1}_{11} =0 \quad \Rightarrow \quad p^{i}_{1i} =0 \label{eq:bip}
 \end{align}
 for all $i=1,2,\dots,d-1$.
\end{thm}

The dual of (\ref{eq:dualNbip}), i.e., $p^{1}_{11} \not=0  \Rightarrow p^{i}_{1i} \not=0, 1\le i \le d-1$, holds for any distance-regular graph \cite[p.178]{BCN}.
By \cite[Corollary 6.2.4]{Dickie},  the graph $\Om$ in Theorem \ref{T:D711} is locally strongly regular:  if for $x\in X$,
 $\Om(x)$ is the set of vertices adjacent to $x$, the induced graph on $\Om(x)$ is strongly regular.

 If $\varGamma$ has $q^1_{11}=0$ (or $\td{q}^1_{11} =0)$, then, by \eqref{eq:dualbip}, it is dual bipartite in case $q^d_{1d} =0$ (or $\td{q}^1_{11} =0)$ and almost dual bipartite otherwise.

   The following result applies when $q^{1}_{11} \ne 0$ and $\tilde{q}^{1}_{11}\ne 0$.

  \begin{lemma} \cite[Lemma 8.2.1]{Dickie}
  \label{L:D821}
Let $\Om$ be distance-regular graph with diameter $d\ge 4$. Suppose
$E_0, E_1,\dots, E_d$ and $\tilde{E}_0, \tilde{E}_1, \dots, \tilde{E}_d$ are $Q$-polynomial structures for $\Om$, with Krein parameters $q^k_{ij}$ and $\tilde{q}^k_{ij}$, respectively. If $q^1_{11}\ne 0$ and $\tilde{q}^{1}_{11}\ne 0$, then  $E_1 = \tilde{E}_d$, $E_d =\tilde{E}_1$ and $d=4$.
\hfill $\Box$
\end{lemma}

By Lemma \ref{L:D821}, $\td{E}_1 = E_4$ and $\td{E}_4 = E_1$.
By Theorem \ref{t:2Q}, the  $Q$-polynomial structures for $\varGamma$ have type III.
The following hold by the $Q$-polynomial property (see also \cite[Theorem 2]{Suzuki98}):
$$
q^{4}_{14} =0 = q^{4}_{34}, \quad q^{4}_{24} \ne 0 \ne q^{4}_{23}.
$$

Pascasio \cite{Pascasio99} showed the following results.
\begin{thm} \cite[Theorem 1.3, Lemma 3.2]{Pascasio99}
\label{T:P13}
Let $\Om$ be a distance regular graph with diameter $d\ge 3$ and eigenvalues $\theta_0 > \theta_1 > \dots > \theta_d$. Let $E$ and $F$ be two primitive idempotents other than $E_0$.
\begin{enumerate}[(i)]
\item
Suppose $\Om$ is tight. Then $E\circ F$ is a scalar multiple of a primitive idempotent $H$ of $\Om$ if and only if $E,F$ are a permutation of $E_1, E_d$. Moreover,  the scalar is $\frac{m_E m_F} {|X| m_H}$.
\item
Suppose $\Om$ is bipartite.  Then $E\circ F$ is a scalar multiple of a primitive of $\Om$ if and only if at least one of $E, F$ is equal to $E_d$.
\item
Suppose $\Om$ is neither bipartite nor tight. Then $E\circ F$ is never a scalar multiple of a primitive of $\Om$.
\hfill $\Box$
\end{enumerate}
\end{thm}

\begin{thm} \cite[Theorem 1.3]{Pascasio01} \label{t:tight}
Let $\Om$ be a $Q$-polynomial distance-regular graph with diameter $d\ge 3,$ intersection numbers $a_i$ and Krein parameters $a_i^*$, $0\le i\le d$. The following are equivalent.
\begin{enumerate}[(i)]
\item
$\varGamma$ is tight.
\item
$\varGamma$ is not bipartite and $a_d =0$.
\item
$\varGamma$ is not bipartite and $a^*_d=0$.
\end{enumerate}
\end{thm}

Since $q^{4}_{14} =0$, $E_1\circ E_4 = |V\varGamma|^{-1} b^*_3 E_3$, where $V\varGamma$ is the vertex set of $\varGamma$. Theorem \ref{T:P13} says that $\theta_4$ is the eigenvalue associated with $E_4$ or $\td{E}_4$. Without loss of generality, we assume that it is $E_4$. If $\varGamma$ is tight, then $E_1$ is associated with $\theta_1$. By \cite[Theorem 1.5]{Pascasio01}, $E_0, E_1, E_2, E_3, E_4$ is the natural ordering of the primitive idempotents. Now we denote the eigenvalues of $\varGamma$ by  $\theta_0 > \theta_1 > \theta_2 > \theta_3 >\theta_4$.

Now we prove Theorem \ref{t:main}  by showing that  one of $a_1^*, \tilde{a}_1^*$ for $\varGamma$ vanishes.  We distinguishing whether $\varGamma$ is bipartite or not.

\subsubsection{$\varGamma$ is bipartite}

Assume that $\varGamma$ is bipartite.  So we have $\theta_4 = - k$ and $m_4 =1$. Imprimitive $Q$-polynomial association schemes have the following characterization.

\begin{thm}(\cite[Theorem 3]{Suzuki98b}, \cite{Cerzo,Tanaka})
\label{t:impQ}
Let $E_0, E_1, \dots, E_d$ be a $Q$-polynomial structure for association  scheme $\sX$.  Suppose that $\sX$ is imprimitive. More precisely, let $T$ be a proper subset of  $\{0,1,\dots, d\}$ with $T\ne \{0\}$ such that the linear span of  $\{E_i\mid i\in T\}$ is closed under the Schur product. In addition, assume $m_1 >2$. Then one of following holds:
\begin{enumerate}[(i)]
\item \label{dualp}
$T=\{0,2,4,\dots\}$ and $a^*_i =0$.
\item \label{dualant}
$T=\{0,d\}$ and $b^*_i = c^*_{d-i}$ for all $i=0,1, \dots, d$ with the possible exception $i=\lfloor d/2\rfloor$.
\end{enumerate}
\end{thm}

An association scheme $\sX$ in case (\ref{dualp}) and (\ref{dualant}) is also called {\em dual bipartite}  and {\em dual antipodal} respectively.

Now back to $\varGamma$.   Let $m_1 = m_{E_1}$.
  Suppose $m_1 >2$. Then Theorem \ref{t:impQ}  applies.  If $q^{1}_{11} >0$, $\varGamma$ can  not be dual bipartite and hence $\varGamma$ is dual antipodal, i.e., case (ii). So  $T=\{0,4\}$. However, $E_4$ is the primary  idempotent for the second $Q$-polynomial structure, which is impossible.

 Suppose $m_1 \le 2$. Since $m_{E_1} < k =3$, we have $\theta_{E_1} = \theta_1$ by \cite[Theorem 4.4.4]{BCN}.  If $m_1 =2$, then by \cite[Theorem 13 (i)]{Koolen13} $k=2$, this contradicts $k>2$. By \cite[Lemma 7]{Koolen13}, it is impossible for  $m_1 =1$; otherwise $m_1 + m_4 =2 < k$.

\subsubsection{$\varGamma$ is not bipartite}

Assume that $\varGamma$ is not bipartite. Since $a^*_4=0$, $\varGamma$ is tight   and  $a_4=0$ by Theorem \ref{t:tight}.  In the literature \cite[p.247]{BCN}, there is an infinite series of feasible formally self-dual intersection arrays
\begin{equation} \label{selfdual}
\{\mu(2\mu +1), (\mu-1)(2\mu +1),\mu^2,\mu; 1,\mu, \mu(\mu-1),\mu(2\mu+1)\}
\end{equation}
This series was ruled out in \cite{Godsil95}.
Had a graph with this array existed, it would be tight with a pair of non-integral eigenvalues,  and would have possessed two $P$-polynomial and two $Q$-polynomial structures. We will show below that there are no distance-regular graphs with intersection array \eqref{selfdual} in an alternative way.

Since $\varGamma$ is not bipartite, it follows from \eqref{eq:dualNbip} that   $a_1a_2a_3 \ne 0$,  $a_4 =0$.
We can deduce from Theorem \ref{T:P13} (i) and Lemma \ref{cosineseq} (iii) that
 \begin{equation}\label{eq:theta14}
 \theta_1 \theta_4 = \theta_0 \theta_3.
 \end{equation}

\begin{thm}\cite[Theorem 8.1.2, Corollary 8.1.4]{BCN}
\label{t:T8.1.2}
Let $\Om$ be a $Q$-polynomial distance-regular graph. Then every $Q$-sequence $(\sg_0, \sg_1, \dots,\sg_d)$ of $\Om$ satisfies the recurrence
\begin{equation}\label{eq:10} 
\sg_{i+1} + \sg_{i-1} = p\sg_i + r \quad (i=1,\dots, d-1),
\end{equation}
for suitable numbers $p$ and $r$.

If $k=\theta_0, \theta_1, \dots, \theta_d$ the $Q$-polynomial structure corresponding to the above $Q$-sequence, then there are constants $r^*, s^*$ such that
\begin{equation} \label{eq:24} 
\left.
\begin{array}{l}
\theta_{\ell +1} + \theta_{\ell -1} = p \theta_{\ell} + r^* \\
\theta_{\ell +1} \theta_{\ell -1}  = \theta_{\ell}^2 - r^* \theta_{\ell} - s^*
\end{array}\right\}
\quad (\ell =1, \dots, d-1).
\end{equation}
\end{thm}

 Now applying Theorem \ref{t:T8.1.2} to the two $Q$-polynomial structures of $\varGamma$, we obtain:
\begin{align}
\theta_2 - p\theta_1 + \theta_0=\theta_3 - p\theta_2 + \theta_1= \theta_4 - p\theta_3 + \theta_2,\label{eq:q1}
\\
\theta_2 - \td{p}\theta_4 + \theta_0=\theta_3 - \td{p}\theta_2 + \theta_4= \theta_1 - \td{p}\theta_3 + \theta_2.\label{eq:q2}
\end{align}
We obtain from \eqref{eq:q1} and \eqref{eq:q2}
$$
p = \frac{\theta_0 - \theta_4}{\theta_1 - \theta_3} ,
\quad
\td{p} = \frac{\theta_0 - \theta_1}{\theta_4 - \theta_3},
\quad
 p - \td{p} = \frac{2(\theta_1 - \theta_4)}{\theta_2 - \theta_3}.
$$
From these equation and $\theta_1\theta_4= k \theta_3$ we find
\begin{align}\label{eq:p}
p = \theta_0/\theta_1,\quad \td{p} = \theta_0/\theta_4, \quad \theta_2 = -\theta_3.
\end{align}
Now substituting these into \eqref{eq:q1} leads to
\begin{equation}\label{eq:theta2}
\theta_1 + \theta_4 = 2\theta_2.
\end{equation}
If we substitute \eqref{eq:p} into \cite[Corollary 8.1.4 (24)]{BCN}, we find
$$
r^* = \theta_2, \quad s^* = \theta_1 \theta_3, \quad
\td{r}^* = \theta_2, \quad \td{s}^* = \theta_4 \theta_2.
$$
Since $\varGamma$ is tight, we have
\begin{equation}
\left(\theta_1 + \frac{k}{a_1 +1}\right)
\left(\theta_4 + \frac{k}{a_1 +1}\right) = -\frac{k a_1 b_1}{(a_1+1)^2}.
\end{equation}
We find from this and Eq. \eqref{eq:theta2} and \eqref{eq:theta14} that
\begin{equation}\label{eq:a1b1}
\theta_2 (a_1 -1) = b_1 +1.
\end{equation}


Now we collect some equations above that are key to the proof to follow:
\begin{align}
\theta_1 \theta_4 &= \theta_0 \theta_3 \label{k1}\\
\theta_1 + \theta_4 &= 2\theta_2 \label{k2}\\
\theta_2 (a_1 -1) &= b_1 +1. \label{k3}
\end{align}
Since $\varGamma$ is tight, a local graph $\varGamma(x)$ is strongly regular with $k$ vertices and valency $a_1$ and non-trivial eigenvalues
 $$\xi = -1 - \frac{b_1}{(\theta_4 +1)},\quad
\tau = - 1- \frac{b_1}{(\theta_1 +1)}.
$$
where $\xi\ge 0$ and $\tau < -1$.

The local graph  $\varGamma(x)$ can not be a conference graph. Otherwise, we have $a_1 = (k-1)/2$ and such a graph has diameter 3 by \cite{Koolen12}. (The intersection array \eqref{selfdual} has the second largest and minimal eigenvalues non-integral, and thus any graph with this array has a conference graph as its local graph and therefore can not exist.)
Therefore, $\sigma$ and $\tau$ are both integers and thus
 $\theta_4, \theta_1$ are both rational numbers. Since they are algebraic integers, $\theta_4, \theta_1$ are  integers.

We find from \eqref{k1}-\eqref{k3} that
$$
-(\theta_1 +1)(\theta_4 +1) =-\theta_1 \theta_4 - \theta_1 - \theta_2 - 1 = (k-2)\theta_2 -1 = (a_1-1) \theta_2 + b_1 \theta_2 -1 = b_1+1 + b_1 \theta_2 -1 = b_1(\theta_2 +1).
$$
From this we can derive
 \begin{equation}
 \dfrac{-b_1^2}{(\theta_1 +1)(\theta_4 +1)} =\dfrac{ b_1}{(\theta_2 +1)}.
 \end{equation}
Since the left hand side is an integer,  $\theta_2 +1$ divides $b_1$ and hence $a_1$ by \eqref{k3}.

Let $(\sg_i)_i$ be the cosine sequence of $\theta_1$. Then $\sg_i > \sg_{i+1}~(0\le i\le 3)$ and $\sg_3=\sg_1\sg_4$. As the sequence $\sg_i$ has one sign change, $\sg_4 <0$ and $\sg_1 >0$. Hence $\sg_3 <0$.
Let $(\td{\sg}_i)_i$ be cosine sequence of $\theta_4$. Then $(-1)^i \td{\sg}_i >0$.  Let $(u_i)_i$ be the cosine sequence of $\theta_3$. Then $u_i = \sg_i\td{\sg}_i$.

 Let $a_1 = \alpha(\theta_2 +1)$,
      $b_1 = \beta(\theta_2 +1)$.
By \eqref{k1} and \eqref{k2} we find $\theta_1, \theta_4 = \theta_2 \pm \sqrt{\theta_2(k + \theta_2)}$.
Now $\theta_2(k+ \theta_2) = \alpha \theta_2(\theta_2  +1)^2$  and hence $k = \alpha(\theta_2+1)^2 - \theta_2$.

By \cite[Proposition 1 (iii), Proposition 2 (iii)]{Cioba},  $\sg_2 \geq 0$ and thus $u_2\ge 0$. We find from this and $\theta_2 = -\theta_3$ that  $k - \theta_2 a_1 \leq \theta_2^2$.
As $k = (a_1 -1)(\theta_2 +1) +1$, we obtain $a_1 \leq \theta_2(\theta_2 +1)$ and thus
$\alpha \leq \theta_2$.

Recall $(\alpha (\theta_2 +1) -1) \theta_2 - 1=(a_1 -1)\theta_2  -1 = b_1 = \beta(\theta_2 +1)$.
We obtain: $\beta= \alpha \theta_2 -1$ and $b_1 = (\alpha \theta_2 -1)(\theta_2 +1)$,
$\xi = -1 - b_1/(\theta_4 +1) = \sqrt{\alpha \theta_2}$ and
$\tau = -\xi$.
This means that $\xi$ divides $a_1 (= k - b_1 -1)$ and hence, as $a_1 = \xi^2 + \alpha$,
$\xi$ divides $\alpha$.
This implies that $\theta_2$ divides $\alpha$ and hence $\theta_2\le \alpha$. We conclude $\alpha = \theta_2 = \xi$. Now,
the local graph $\varGamma(x)$ has the following parameters:
$$
k = \theta_2^2(\theta_2 +2), \quad
a_1 = \theta_2 (\theta_2 +1), \quad
\theta_1 = \theta_2(\theta_2 +2), \quad
 \theta_4 = -\theta_2^2.
$$
We see that
$$\theta_1 = \frac{a_1 + \sqrt{a_1^2 + 4k}}{2}.
$$
Now $\varGamma$ is antipodal by the following result.

\begin{lemma} \cite[Proposition 3.5]{Koolen11} \label{l:prop3.5}
Let $\Om$ be a distance-regular graph with $d$ at least three and distinct eigenvalues
$k=\theta_0 > \theta_1 > \dots >\theta_d$. Then $\theta_1 = ({a_1 + \sqrt{a_1^2 + 4k}})/{2}$ if and only if  one of the following holds:
\begin{enumerate}[(i)]
\item
$d = 3$ and $\Om$ is a Shilla distance-regular graph;
\item
$ d = 4$ and $\Om$ is an antipodal distance-regular graph.
\end{enumerate}
\end{lemma}

Since $\varGamma$ is antipodal, it is dual bipartite. So $a_1^*=0$, or $\tilde{a}^*_1=0$. This completes the proof of Theorem \ref{t:main}.

\subsection{Proof of Theorem \ref{mainD=4}}

 By Lemma \ref{T:D512} and the succeeding remark, the parameters $a_i^*$ of $\varGamma$ can be divided into the following four cases:
\begin{enumerate}[(i)]
\item
If $a_1^* =a_4^*=0$, then $\varGamma$ is $H(4,2)$ or a Hadamard graph by \cite{Dickie96}.
\item
If $a_1^* = 0\ne a_4^*$, then $\varGamma$ is  $\frac{1}{2}H(9,2)$, or $\td{H}(9,2)$ by \cite[Theorem 3.1.4]{Dickie}.

\item
Suppose $a_1^* \ne 0\ne a_4^*$. By \cite[Theorem 2]{Suzuki98},  the $Q$-polynomial structures $(E_i)_i$ is almost dual antipodal and thus by \cite[Lemma 3.1.3]{Dickie}, the other $Q$-polynomial structure $(\tilde{E}_i)_i$ is  almost dual bipartite. This case is implied by the previous case by treating  $(\tilde{E}_i)_i$, which has  $\tilde{a}_1^*= 0 \ne \tilde{a}_4^*$.
\item
If $a_1^* \ne 0 = a_4^*$, then $\tilde{a}_1^*=0$ by Theorem \ref{t:main} and Cases (i),(ii) apply.
\end{enumerate}
So the proof of Theorem \ref{mainD=4} is completed.

 \begin{rmk}
 The following graphs are  twice $Q$-polynomial distance-regular graphs  of diameter $3$: each graph in Theorem \ref{t:Dickie} (ii)-(iv) with $d=3$, and a  distance-regular graph intersection array
 $
 \{k, k - a_1 -1, 1; 1,k -a_1 -1, k\}.
 $
 If $a_1=0$,  this array is uniquely realized by the complement of $K_{k+1}\times K_2$; any distance-regular graph with this intersection array is called a {\em Taylor graph} if $a_1 >0$; see \cite[p.13]{BCN}.
 \end{rmk}

\section*{Acknowledgments}

J. Koolen thanks the 100 Talents Program of the Chinese Academy of Sciences for support.
J. Ma acknowledges  support from the Natural Science Foundation of Hebei province (A2012205079) and Science Foundation of  Hebei Normal University (L2011B02).

 
\end{document}